\documentclass[11pt]{article}
\usepackage{amsmath,amssymb,amsthm,url}
\usepackage{color}
\usepackage[margin=1.25in,letterpaper]{geometry}
\usepackage{graphicx}

\theoremstyle{definition}

\newcommand{\FF}{\mathbb{F}}

\DeclareMathOperator{\Aut}{Aut}
\DeclareMathOperator{\PG}{PG}
\DeclareMathOperator{\AG}{AG}

\title{On the $486$-vertex distance-regular graphs\\ of Koolen--Riebeek and Soicher}
\author{Robert F.\ Bailey\footnote{School of Science and the Environment (Mathematics), Grenfell Campus, Memorial University of Newfoundland, Corner Brook, NL A2H~6P9, Canada. E-mail: \texttt{rbailey@grenfell.mun.ca}} \and Daniel R.\ Hawtin\footnote{Department of Mathematics, University of Rijeka, Radmile Matej\v{c}i\'c 2, 51000 Rijeka, Croatia. E-mail: \texttt{dhawtin@math.uniri.hr}}}

\begin{document}
\maketitle

\begin{abstract}
This paper considers three imprimitive distance-regular graphs with $486$ vertices and diameter~$4$: the Koolen--Riebeek graph (which is bipartite), the Soicher graph (which is antipodal), and the incidence graph of a symmetric transversal design obtained from the affine geometry $\AG(5,3)$ (which is both).  It is shown that each of these is preserved by the same rank-$9$ action of the group $3^5:(2\times M_{10})$, and the connection is explained using the ternary Golay code.\\



\end{abstract}

\section{Introduction} \label{sec:intro}

A (finite, simple, connected) graph with diameter~$d$ is {\em distance-regular} if, for all $i$ with $0\leq i\leq d$ and any vertices $u,w$ at distance $i$, the number of neighbours of $w$ at distances $i-1$, $i$ and $i+1$ from $u$ depends only on $i$, and not on the choices of $u$ and $w$.  These numbers are denoted by $c_i$, $a_i$ and $b_i$ respectively, and are known as the {\em parameters} of the graph.  It is easy to see that $c_0$, $b_d$ are undefined, $a_0=0$, $c_1=1$, and that the graph must be regular with degree $b_0=k$.  Consequently $c_i+a_i+b_i=k$ for $1\leq i\leq d-1$ and $c_d+a_d=k$, so the parameters $a_i$ are determined by the others.  We put the parameters into an array, called the {\em intersection array} of the graph,
\[ \{ k,b_1,\ldots,b_{d-1}; \, 1,c_2,\ldots,c_d \}. \]
Two important class of distance-regular graphs are the {\em distance-transitive} graphs, where a group of automorphisms of the graph acts transitively on pairs of vertices at each distance, and the {\em strongly regular graphs}, which are the distance-regular graphs of diameter~$2$. In this case, the parameters are usually given as $(n,k,\lambda,\mu)$, where $n$ is the number of vertices, $\lambda=a_1$ and $\mu=c_2$.

A distance-regular graph is called {\em primitive} if its distance-$i$ graphs are connected for $1\leq i\leq d$, and {\em imprimitive} otherwise; it is well-known that the only possibilities for a distance-regular graph to be imprimitive are either if it is bipartite, or it is {\em antipodal}, which occurs when then distance-$d$ graphs consist of disjoint cliques (although both possibilities may occur in the same graph).  An imprimitive distance-regular graph can be reduced to a primitive one by the operations of {\em halving} and/or {\em folding}: the halved graphs are the connected components of the distance-$2$ graph of a bipartite graph; the folded graph is the quotient graph obtained on the antipodal classes.  An antipodal graph may be referred to as an {\em antipodal cover} of its folded graph.  For further background, we refer to the book of Brouwer, Cohen and Neumaier~\cite{BCN} and the survey of van Dam, Koolen and Tanaka~\cite{vDKT}.

\section{The graphs of interest} \label{sec:graphs}

In this paper, we are primarily concerned with the following graphs.

For a linear code $\mathcal{C}$ of dimension $k$ in $\FF_q^n$, the {\em coset graph} has as its vertices the $q^{n-k}$ cosets of $\mathcal{C}$, with two cosets being adjacent if they contain representatives at Hamming distance~$1$.  The {\em ternary Golay code} $\mathcal{G}$ is a $6$-dimensional linear code in $\FF_3^{11}$, which is the unique such perfect code with minimum distance~$5$ (see~\cite{MacWS77}); its coset graph is the {\em Berlekamp--van Lint--Seidel graph}, obtained in~\cite{BvLS73}, which is a strongly regular graph $\Gamma$ with parameters $(243,22,1,2)$.  The complement of $\Gamma$ is strongly regular with parameters $(243,220,199,200)$.  The automorphism group of $\Gamma$ has the form $3^5:(2\times M_{11})$, with rank~$3$ (i.e.\ there are three orbits on ordered pairs of vertices); consequently $\Gamma$ is distance-transitive.

The {\em Koolen--Riebeek graph} is a bipartite distance-regular graph $\Delta$ with $486$ vertices, diameter~$4$ and intersection array
\[ \{45,44,36,5;\, 1,9,40,45\}. \]
Described in the 1998 paper of Brouwer, Koolen and Riebeek~\cite{BKR98}, its halved graphs are isomorphic to the complement of the Berlekamp--van Lint--Seidel graph.  Its automorphism group has the form $3^5:(2\times M_{10})$, but as its rank is~$9$, the graph is not distance-transitive.  In~\cite{BKR98}, a description of the graph is given as the incidence graph of an incidence structure whose points are the cosets of the ternary Golay code $\mathcal{G}$ and whose blocks are a particular collection of $45$-cocliques in the Berlekamp--van Lint--Seidel graph $\Gamma$.

The {\em (second) Soicher graph} $\Upsilon$ is one of three distance-regular graphs given in Soicher's 1993 paper~\cite{Soicher93}.  It also has $486$ vertices, and has intersection array
\[ \{56,45,16,1;\, 1,8,45,56\}. \]
It is an antipodal graph; its folded graph is the unique strongly regular graph with parameters $(162,56,10,24)$, the second subconstituent of the {\em McLaughlin graph}.  Soicher constructed the graph computationally, starting from the Suzuki simple group.  Its automorphism group has the form $3\cdot U_4(3):2^2$ and has rank~$5$, so $\Upsilon$ is distance-transitive.  In an unpublished manuscript, Brouwer showed that $\Upsilon$ is the only distance-regular graph with these parameters.  (Brouwer's proof can be found in~\cite{BCN-corr-add}.)

The bipartite, antipodal distance-regular graphs of diameter~$4$ are precisely the incidence graphs of {\em symmetric transversal designs} (also known as {\em symmetric nets}: see~\cite[{\S}1.7]{BCN} or~\cite{incidence} for details).  From the affine geometry $\AG(n,q)$, one may obtain a symmetric transversal design by choosing two arbitrary points, then deleting each parallel class of $(n-1)$-flats containing a flat through those two points (see~\cite[Proposition 7.18]{BJL} for details).  The incidence graph of this design has $2q^n$ vertices and intersection array
\[ \{ q^{n-1},q^{n-1}-1,q^{n-1}-q^{n-2},1;\, 1,q^{n-2},q^{n-1}-1,q^{n-1} \}. \]
Such a design has automorphism group with index $(q^n-1)/(q-1)$ in $\mathrm{A\Gamma L}(n,q)$ (see~\cite{Jungnickel84}); this group is an index-$2$ subgroup of the automorphism group of the incidence graph.
These graphs are in fact distance-transitive: see Ivanov {\em et al.}~\cite{ILPP97} for an alternative construction from the projective space $\PG(n,q)$.  
However, the graphs are not determined by their parameters: for example, there are exactly four graphs with the parameters of that arising from $\AG(3,3)$~\cite{Mavron2000}, while more generally there are vast numbers of non-isomorphic graphs with these parameters (cf.~\cite{Jungnickel84}).
In this paper, we are interested in the design and graph obtained from the $4$-flats in $\AG(5,3)$, where we have a distance-transitive graph $\Sigma$ with $486$ vertices and intersection array
\[ \{81,80,54,1;\, 1,27,80,81\}. \]

The main result of this paper is that the graphs $\Delta$, $\Upsilon$ and $\Sigma$ can each be constructed from the same rank-$9$ action of the group $3^5:(2\times M_{10})$.  At first, this was observed computationally, but we were later able to obtain a theoretical explanation of this observation.

\section{Computer construction of the graphs} \label{sec:computer}
The first author, along with his students, has been developing an online catalogue of distance-regular graphs~\cite{distanceregular.org}, using the {\sf GRAPE} package~\cite{GRAPE} for the {\sf GAP} computer algebra system~\cite{GAP}.  Typically, these graphs are obtained by providing a group of automorphisms, then using the {\small \tt EdgeOrbitsGraph} function in {\sf GRAPE}.
To obtain constructions of the graphs $\Gamma$, $\Delta$ and $\Upsilon$ to include in the catalogue, the following approach was used.

First, to construct the Koolen--Riebeek graph $\Delta$, we use the fact that $\Aut(\Gamma)\cong 3^5:(2\times M_{11})$ is a primitive permutation group of rank $3$.  In the {\sf GAP} library of primitive groups, there are in fact two groups of degree $243$ with this structure.  Using {\small \tt EdgeOrbitsGraph}, we find that one is the automorphism group of $\Gamma$, and the other preserves a strongly regular graph with parameters $(243,110,37,60)$ (the {\em Delsarte graph}, which also arises from the ternary Golay code: see~\cite{BvL84}).  Let $G=\Aut(\Gamma)$.  Now, $H=3^5:(2\times M_{10})$ may be obtained as a maximal subgroup of index $11$ in $G$ in {\sf GAP}.  Next, since we want a transitive action of $H$ of degree~$486$, we find the conjugacy classes of subgroups of $H$ with index~$486$, and examine the actions of $H$ on right cosets.  We find that there are four such actions of degree~$486$, but only one has rank~$9$.  Generators for this rank-$9$ action of $H$ are given in the Appendix.

The {\small \tt VertexTransitiveDRGs} function in {\sf GRAPE}, using the technique of collapsed adjacency matrices described in~\cite{PS97}, determines the intersection arrays of distance-regular graphs with a given vertex-transitive group of automorphisms.  The intersection arrays, and corresponding graphs, for the rank-$9$ action of $H$ obtained are as follows:
\renewcommand{\arraystretch}{1.33}
\[ \begin{array}{ll}
\textnormal{Intersection array} & \textnormal{Graph} \\ \hline
\{485;\, 1\}                & K_{486} \\
\{243,242;\, 1,243\}        & K_{243,243} \\
\{483,2;\, 1,483\}          & K_{3^{162}} \\
\{45,44,36,5;\, 1,9,40,45\} & \Delta \\
\{56,45,16,1;\, 1,8,45,56\} & (\ast)\\
\{81,80,54,1;\, 1,27,80,81\} & (\dagger)
\end{array} \]%
\renewcommand{\arraystretch}{1.0}%
The edge set of each graph is obtained as a union of orbitals of $H$, i.e.\ orbits on pairs of vertices.  Now, there is a single orbital of $H$ where the corresponding suborbit has length $45$; this orbital must give the edges of $\Delta$, so we may use {\small \tt EdgeOrbitsGraph} to construct the graph in {\sf GRAPE}.  We can also use {\sf GRAPE} to confirm that the full automorphism group of $\Delta$ is precisely $H$.

However, we notice that the intersection array $(\ast)$ is precisely that of the second Soicher graph $\Upsilon$.  Since $\Upsilon$ is the unique distance-regular graph with this intersection array, then that must be the graph we have obtained here.  In addition, we observe that the intersection array $(\dagger)$ is precisely that of the graph $\Sigma$ obtained from the $4$-flats of $\AG(5,3)$.  It is straightforward to construct $\Sigma$ in {\sf GRAPE} (beginning with the {\sf DESIGN} package~\cite{DESIGN} to obtain the $4$-flats in $\AG(5,3)$) and to confirm that the graph obtained from $H$ really is isomorphic to $\Sigma$.  To the best of our knowledge, the existence of either $\Upsilon$ or $\Sigma$ as a union of orbitals of $\Aut(\Delta)$ was not previously known.

\subsection{Orbit diagrams}
\begin{figure}[!hbt]
\centering
 \includegraphics[scale=0.8]{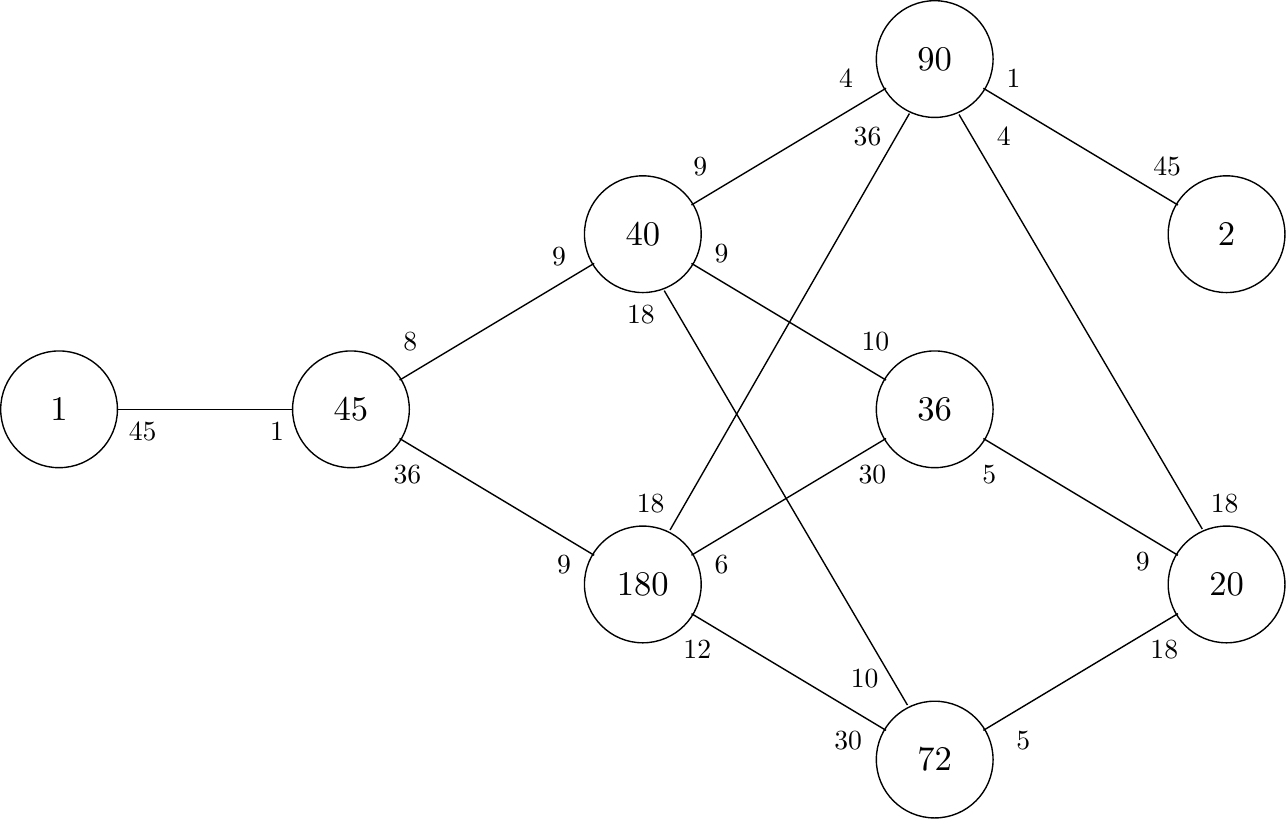}
 \caption{Orbit diagram for the Koolen-Riebeek graph $\Delta$ relative to $3^5:(2\times M_{10})$.}
 \label{orbitdiagKR}
 
 \vspace*{\floatsep}
 
 \includegraphics[scale=0.8]{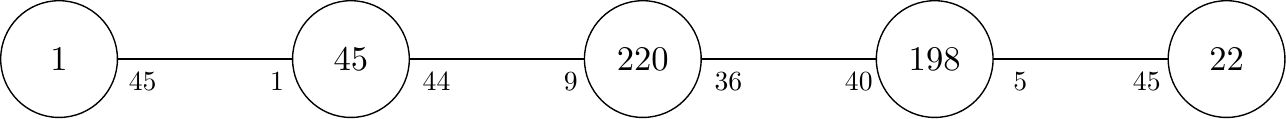}
 \caption{Distance distribution diagram for the Koolen-Riebeek graph $\Delta$.}
 \label{distdiagKR}
 
 \vspace*{\floatsep}
 
 \includegraphics[scale=0.8]{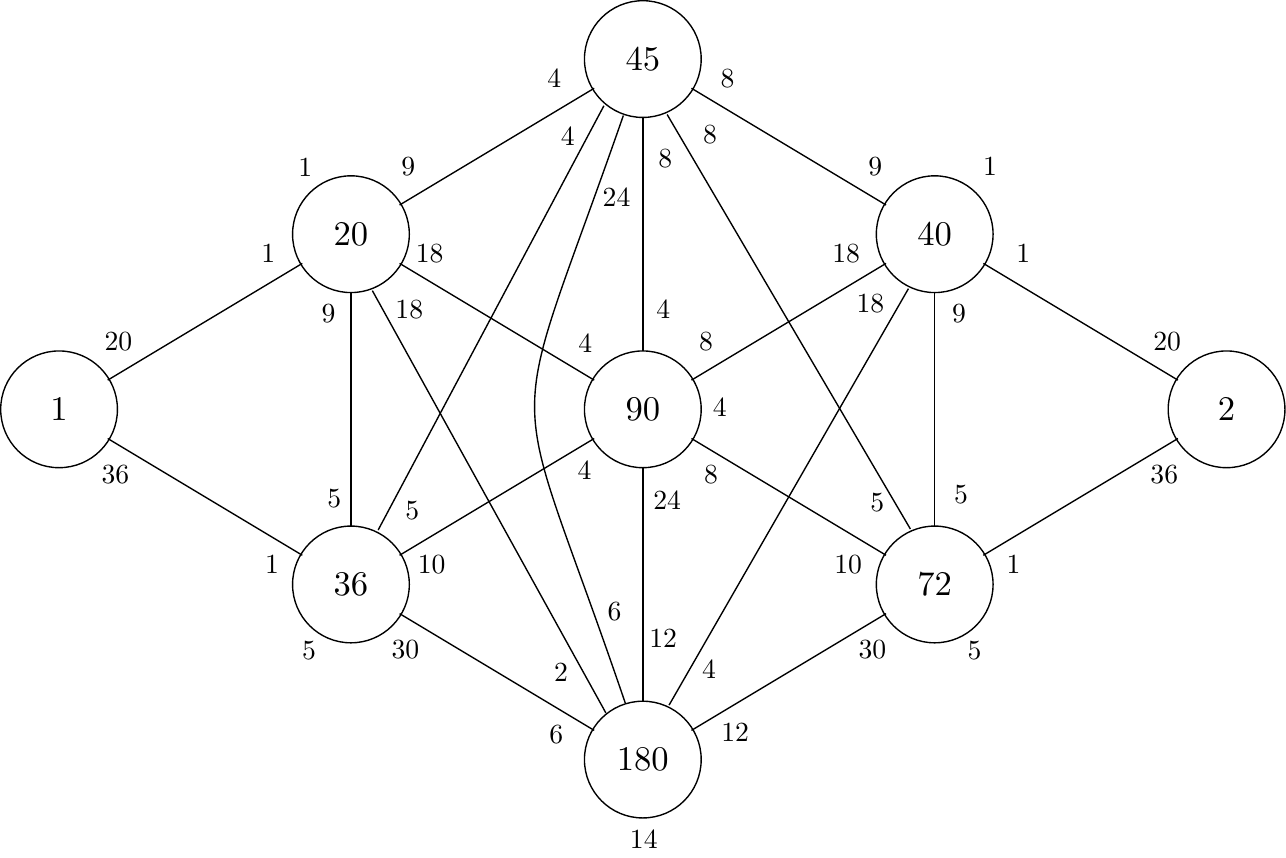}
 \caption{Orbit diagram for the Soicher graph $\Upsilon$ relative to $3^5:(2\times M_{10})$.}
 \label{orbitdiagSoicher}
 
 \vspace*{\floatsep}
 
 \includegraphics[scale=0.8]{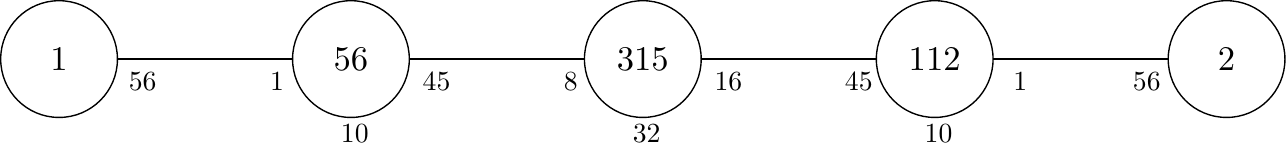}
 \caption{Distance distribution diagram for the Soicher graph $\Upsilon$.}
 \label{distdiagSoicher}

\end{figure}

\begin{figure}[!hbt]
\centering

 \includegraphics[scale=0.8]{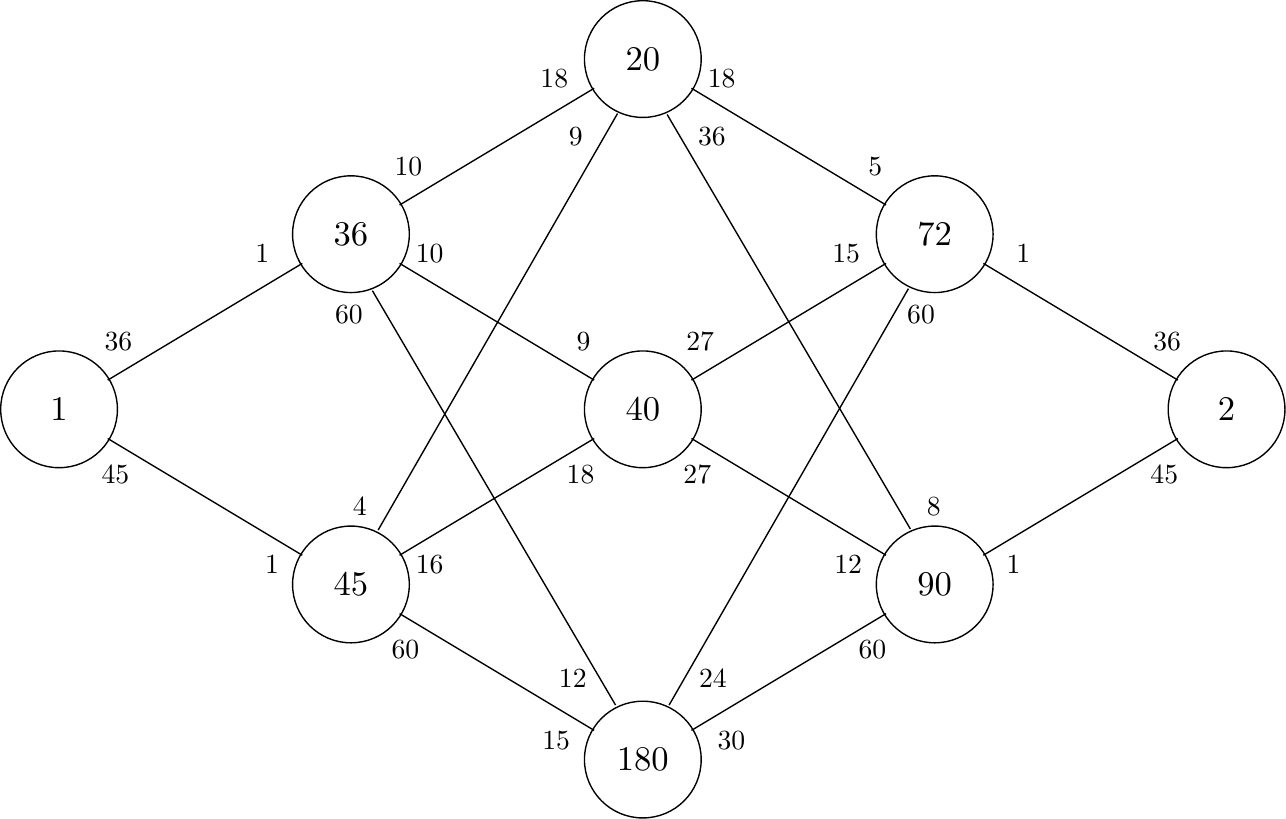}
 \caption{Orbit diagram for the incidence graph $\Sigma$ relative to $3^5:(2\times M_{10})$.}
 \label{orbitdiagSTD}
 
 \vspace*{\floatsep}
 
 \includegraphics[scale=0.8]{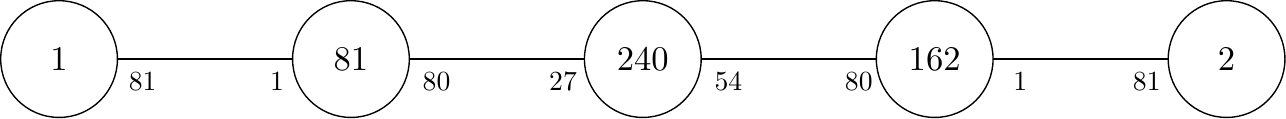}
 \caption{Distance distribution diagram for the incidence graph $\Sigma$.}
 \label{distdiagSTD}
\end{figure}

The group $H\cong 3^5:(2\times M_{10})$ has $9$ suborbits, of lengths $1,2,20,36,40,45,72,90,180$.  To visualize how the three graphs $\Delta$, $\Upsilon$ and $\Sigma$ arise from the suborbits, we can examine their orbit diagrams (Figures~\ref{orbitdiagKR}, \ref{orbitdiagSoicher} and~\ref{orbitdiagSTD}).  Each node corresponds to a suborbit, labelled by its size; the edge labels correspond to how many neighbours a vertex in a given suborbit has in the adjacent one.  These edge labels are precisely the entries of the collapsed adjacency matrix of the graph relative to the group of automorphisms $H$ (see~\cite[{\S}2.3]{PS97}), which can be calculated in {\sf GRAPE}.  Also, in Figures~\ref{distdiagKR}, \ref{distdiagSoicher} and~\ref{distdiagSTD} we give the distance distribution diagrams for each of the three graphs.  (The orbit diagram of $\Delta$ is taken from~\cite{BKR98}.)

%\newpage
\section{Explaining the connection} \label{sec:connection}

The ternary Golay code $\mathcal{G}$ is a $6$-dimensional linear code in $\FF_3^{11}$ with minimum distance~$5$.  Following~\cite{BKR98}, it can be obtained from a generator matrix whose rows are the $11$ cyclic permutations of $(-+-+++---+-)$ (where $+,-$ denote the non-zero elements of $\FF_3$).  Let $V=\FF_3^{11}$ and consider the quotient space $W=V/\mathcal{G}$, whose elements are the cosets of $\mathcal{G}$.  Clearly, the affine space $\AG(W)$ is isomorphic to $\AG(5,3)$; the flats of $\AG(W)$ are the cosets of subspaces of $V$ containing $\mathcal{G}$.  Now, let $\mathcal{F}$ denote the set of $4$-flats of $\AG(W)$ which, when viewed as $10$-flats in $\AG(V)$, are the cosets of $10$-spaces in $V$ containing $\mathcal{G}$ but not $\mathcal{G}+e_0$ (where $e_0,\ldots,e_{10}$ denote the standard basis vectors for $V$).  The number of such subspaces is 
\[ \begin{bmatrix} 11-6 \\ 10-6 \end{bmatrix}_3 - \begin{bmatrix} 11-7 \\ 10-7 \end{bmatrix}_3 = \begin{bmatrix} 5\\4 \end{bmatrix}_3 - \begin{bmatrix} 4\\3 \end{bmatrix}_3 = \frac{3^5-1}{3-1} - \frac{3^4-1}{3-1} = 3^4 = 81, \]
so there are $3^5=243$ flats.  By construction, the incidence graph of the points of $\AG(W)$ (i.e.\ cosets of $\mathcal{G}$) versus the collection $\mathcal{F}$ of flats is precisely the graph $\Sigma$ defined above; in this construction, the neighbours of $\mathcal{G}$ are the $81$ subspaces of $V$ containing $\mathcal{G}$ but not $\mathcal{G}+e_0$.

Now consider this collection of $81$ subspaces; each subspace can be viewed as a $10$-dimensional code in $\FF_3^{11}$ with minimum distance~$1$.  Using {\sf Magma}~\cite{Magma}, one can obtain the weight distributions of these codes, where there are two classes.  There are $45$ subspaces with weight distribution
\[ 0^1\,1^{10}\,2^{70}\,3^{420}\,4^{1770}\,5^{4992}\,6^{9822}\,7^{13980}\,8^{14160}\,9^{9440}\,10^{3680}\,11^{704} \]
which we will refer to as {\em Type~I} subspaces, while there are $36$ subspaces with weight distribution
\[ 0^1\,1^4\,2^{76}\,3^{456}\,4^{1716}\,5^{4956}\,6^{9912}\,7^{13944}\,8^{14214}\,9^{9314}\,10^{3776}\,11^{680} \]
which we will refer to as {\em Type~II} subspaces.  Let $\mathcal{T}_1$ and $\mathcal{T}_2$ denote the sets of all Type~I and Type~II subspaces respectively.  Using {\sf Magma} again, we can calculate the setwise stabilizers of $\mathcal{T}_1$ and $\mathcal{T}_2$ in $\Aut(\Sigma)$: each of these is $2\times M_{10}$.
Next, construct the incidence graph of the cosets of $\mathcal{G}$ versus the $10$-flats $\mathcal{F}$, where $\mathcal{G}$ is incident with the Type~I subspaces, and a coset $\mathcal{G}+x$ is incident with the set of translates $\{U+x\, : \, U\in\mathcal{T}_1\}$.
Another computer calculation (this time in {\sf GAP}) verifies that this is the Koolen--Riebeek graph $\Delta$.  
By construction, the group $3^5:(2\times M_{10})$ acts as automorphisms of this graph.

The stabilizer of $\mathcal{G}$ in $\Aut(\Delta)$ is $M_{10}$, which has five orbits on the cosets of $\mathcal{G}$ of lengths $1,2,20,40,180$.  Since $\mathcal{G}$ is a perfect code with minimum distance~$5$, each coset is given by a unique representative of weight $0$, $1$ or $2$.  Another computation in {\sf Magma} shows that the coset representatives of each orbit can be described as follows.
\renewcommand{\arraystretch}{1.33}
\[ \begin{array}{ll}
\textnormal{Orbit size} & \textnormal{Coset representative} \\ \hline
  1 & \bf{0} \\
  2 & \pm e_0 \\
 20 & \pm e_i~\textnormal{($i\neq 0$)} \\
 40 & \pm e_0 \pm e_i~\textnormal{($i\neq 0$)} \\
180 & \pm e_i \pm e_j~\textnormal{($0<i<j$)}
\end{array} \]%
\renewcommand{\arraystretch}{1.0}%
This action of $M_{10}$ has four orbits on the set $\mathcal{F}$ of $10$-flats of lengths $36,45,72,90$: the orbits of lengths $45$ and $36$ are $\mathcal{T}_1$ and $\mathcal{T}_2$, while the orbits of lengths $90$ and $72$ are formed of the cosets of the $10$-subspaces in $\mathcal{T}_1$ and $\mathcal{T}_2$ respectively.  We can therefore obtain a valency-$56$ graph on the same vertex set as $\Sigma$ and $\Delta$, by putting $\mathcal{G}$ adjacent to the $36$ Type~II subspaces and the $20$ weight-$1$ cosets $\mathcal{G}\pm e_i$ (for $i\neq 0$), and then applying the action of $\Aut(\Delta)\cong 3^5:(2\times M_{10})$.  But this is precisely the construction of the Soicher graph $\Upsilon$ given above.

\section{An induced subgraph} \label{sec:induced}
\begin{figure}[!hbtp]
\centering
  \includegraphics[scale=0.8]{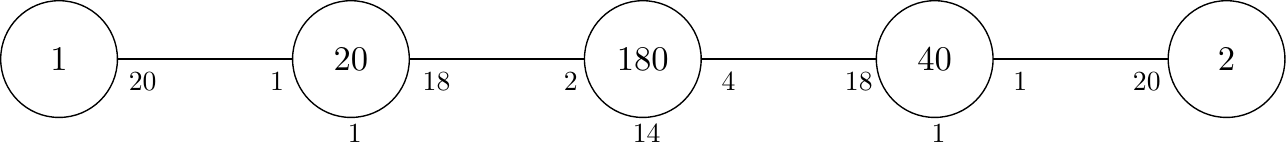}
 \caption{Distance distribution diagram for the induced subgraph $\Lambda$.}
 \label{distdiagInduced}
\end{figure}

The construction of the Soicher graph~$\Upsilon$ from the ternary Golay code divides the vertices into two classes, i.e.\ the bipartition of $\Sigma$ or $\Delta$, formed of the cosets of $\mathcal{G}$ and the collection of $10$-flats.  It seems natural to consider the induced subgraphs of $\Upsilon$ on each of these classes; without loss of generality we consider that on the cosets of $\mathcal{G}$.  By our construction, we can see that $\mathcal{G}$ will be adjacent to the $20$ weight-$1$ cosets $\mathcal{G}\pm e_i$ (for $i\neq 0$).  From the orbit diagram in Figure~\ref{orbitdiagSoicher}, by considering just those suborbits on cosets, we obtain a distance-transitive graph $\Lambda$ with $243$ vertices and intersection array
\[ \{20,18,4,1;\, 1,2,18,20\}. \]
Its distance distribution diagram, which is also its orbit diagram relative to $3^5:(2\times M_{10})$, is given in Figure~\ref{distdiagInduced}.

This graph is also known: it is labelled (A17) in~\cite[{\S}11.3H]{BCN}, where it is shown to be isomorphic to the coset graph of the shortened ternary Golay code (i.e.\ the code obtained from $\mathcal{G}$ by deleting a single co-ordinate from each codeword, and then keeping only those codewords that had a zero in that position).   Furthermore, it is known to be an antipodal $3$-cover of the coset graph of the truncated ternary Golay code (the {\em Brouwer--Haemers graph}~\cite{BH92}, the unique strongly regular graph with parameters $(81,20,1,6)$).  The fact that $\Aut(\Lambda)\cong 3^5:(2\times M_{10})$ is mentioned in~\cite{BCN-corr-add}.  However, we believe that the observation that $\Lambda$ is an induced subgraph of $\Upsilon$ is new.

\section*{Appendix: Generators of $3^5:(2\times M_{10})$}
\label{app:gens}
The generators $a,b,c$ of $H\cong\Aut(\Delta)\cong 3^5:(2\times M_{10})$ were obtained from the {\sf GAP} output of the computations described in Section~\ref{sec:computer}.  A suitable permutation in $\mathrm{Sym}(486)$ was used to obtain a conjugate labelled so that the bipartite halves of $\Delta$ were labelled $1,\ldots,243$ and $244,\ldots,486$.\\

\noindent {\footnotesize \tt %
%%%%%%%%%%%%%%%%%%%%%%%%%%%%%%%%%%%%%%%%%%%%%%%%%%%%%%%%%%%%%%%%%%%%%%%%%%%%%%%%%%%%%%%
a := (1,319,43,476)(2,414,236,469)(3,266,56,295)(4,321,104,290)(5,444)(6,285,237,472)\\
(7,308,179,247)(8,244)(9,447,75,263)(10,468,54,440)(11,389,183,471)(12,293,232,335)\\
(13,314,152,362)(14,364,73,442)(15,423,37,300)(16,382,105,390)(17,478,74,272)\\
(18,342,234,286)(19,407,197,325)(20,316,131,462)(21,301,30,401)(22,402,153,311)(23,483)\\
(24,460,70,434)(25,279,186,367)(26,484,133,404)(27,296,71,410)(28,278,150,357)\\
(29,431,137,466)(31,305,94,317)(32,482,188,416)(33,450)(34,351,72,264)(35,331,127,262)\\
(36,388,235,273)(38,258,191,359)(39,454,223,421)(40,306,201,312)(41,288)\\
(42,370,101,291)(44,328,109,344)(45,458,203,354)(46,256,126,309)(47,313,130,245)\\
(48,349,114,398)(49,396,142,343)(50,393,118,365)(51,441)(52,352,212,403)\\
(53,409,220,456)(55,340,79,449)(57,383,209,381)(58,271,176,392)(59,249)(60,400,83,385)\\
(61,281,178,347)(62,439,87,345)(63,485,80,283)(64,455,69,486)(65,337,86,248)\\
(66,453,81,412)(67,332,85,448)(68,424,78,341)(76,356)(77,324,82,386)(84,399)\\
(88,473,228,333)(89,457,140,427)(90,429,219,425)(91,479)(92,481)(93,463,240,289)\\
(95,430)(96,397,117,438)(97,371,120,467)(98,252,192,446)(99,284,151,254)\\
(100,445,141,426)(102,451,171,339)(103,355,204,375)(106,299,185,287)(107,475,166,251)\\
(108,395)(110,373,184,257)(111,276,144,474)(112,436)(113,477,156,282)(115,372,149,470)\\
(116,255,157,420)(119,422,205,294)(121,277,138,260)(122,394,210,361)(123,298,214,269)\\
(124,418,241,384)(125,419,216,459)(128,330,158,350)(129,452,155,377)(132,480)\\
(134,250,159,310)(135,334,154,417)(136,379)(139,366,217,391)(143,443,187,297)\\
(145,315,162,360)(146,405,224,280)(147,265,242,465)(148,336)(160,261,172,411)\\
(161,387,200,348)(163,303,226,376)(164,464,193,268)(165,378,169,274)(167,307,173,374)\\
(168,368)(170,320,195,318)(174,275,215,432)(175,363)(177,304)(180,246)(181,353)\\
(182,329)(189,435,221,253)(190,408)(194,323,231,437)(196,413,206,428)(198,259,230,380)\\
(199,267,222,322)(202,270)(207,326,239,369)(208,415,213,358)(211,433,238,461)(218,327)\\
(225,406,227,346)(229,292)(233,338,243,302)}\\
	
\noindent {\footnotesize \tt %
%%%%%%%%%%%%%%%%%%%%%%%%%%%%%%%%%%%%%%%%%%%%%%%%%%%%%%%%%%%%%%%%%%%%%%%%%%%%%%%%%%%%%%%
b := (1,90,89,76,36,221)(2,192,50,22,166,220,39,8,176,140,243,45)\\
(3,77,72,231,73,23,177,67,79,164,9,86)(4,209,181,28,161,201,189,163,172,121,185,110)\\
(5,218,56,124,169,214,115,206,26,135,122,180)\\
(6,88,118,183,165,200,188,61,175,239,184,141)\\
(7,179,34,198,25,70,19,157,85,146,91,81)\\
(10,132,87,114,158,139,194,240,100,131,224,120)\\
(11,174,187,13,167,35,191,53,29,155,210,147)\\
(12,102,62,236,31,204,186,94,117,207,195,213)\\
(14,68,143)(15,75,126,47,27,229,112,33,162,197,60,153)(16,228,46,216,168,234)\\
(17,152,190,78,134,235,82,225,32,69,238,149)\\
(18,101,127,104,160,123,37,103,171,92,52,125)\\
(20,106,211,107,49,219,193,96,51,95,108,156)\\
(21,137,40,241,54,144,116,59,173,128,48,119)\\
(24,105,43,64,170,83,80,98,65,84,242,63)(30,199,148,41,111,208,130,142,178,113,203,217)\\
(38,71,150,138,55,99,223,66,227,129,74,93)(42,205,230,44,202,136,97,215,154,196,57,145)\\
(58,151,222,212,232,109)(133,159,182,233,237,226)\\
(244,280,428,469,314,399,266,288,348,463,423,410)\\
(245,352,293,437,319,335,258,474,342,462,447,429)\\
(246,458,351,289,442,409,251,329,296,369,424,347)\\
(247,486,387,323,426,324,250,385,299,338,327,453)\\
(248,441,386,422,466,328,272,470,356,397,465,403)\\
(249,357,456,290,406,476,268,311,344,461,416,455)\\
(252,365,480,326,355,446,262,303,313,433,322,298)\\
(253,445,436,316,332,277,255,432,330,484,337,282)\\
(254,377,448,302,384,413,256,417,341,297,380,333)\\
(257,259,312,370,427,378,274,286,395,307,350,393)(260,483,390,362,345,310)\\
(261,485,301,460,354,419,287,400,388,450,439,349)\\
(263,346,421,353,477,304,283,367,372,438,408,440)(264,359,340,294,305,275)\\
(265,478,482,414,473,334,269,364,402,321,411,479)\\
(267,375,300,368,394,278,271,383,292,291,366,279)(270,435,454,374,459,318)\\
(273,407,343,404,464,392,284,401,376,443,317,467)\\
(276,471,405,449,398,425,281,360,472,434,430,418)(285,339,412,457,389,336)\\
(295,382,306,444,325,371,320,468,415,308,309,381)\\
(315,379,373,431,475,361,481,420,396,331,451,363)(358,391,452)}\\
		
\noindent {\footnotesize \tt %
%%%%%%%%%%%%%%%%%%%%%%%%%%%%%%%%%%%%%%%%%%%%%%%%%%%%%%%%%%%%%%%%%%%%%%%%%%%%%%%%%%%%%%%
c := (1,110,193)(2,15,142)(3,94,91)(4,235,8)(5,97,173)(6,116,128)(7,189,10)(9,225,107)\\
(11,65,221)(12,109,49)(13,17,164)(14,206,99)(16,195,209)(18,92,188)(19,79,236)\\
(20,117,151)(21,169,57)(22,172,152)(23,147,32)(24,76,167)(25,185,158)(26,230,54)\\
(27,217,166)(28,108,36)(29,170,90)(30,224,219)(31,212,51)(33,234,242)(34,73,213)\\
(35,131,208)(37,122,205)(38,59,145)(39,238,64)(40,144,175)(41,176,162)(42,74,125)\\
(43,153,216)(44,55,119)(45,161,69)(46,186,121)(47,228,80)(48,137,165)(50,82,63)\\
(52,56,136)(53,134,77)(58,229,191)(60,86,198)(61,237,138)(62,201,168)(66,141,182)\\
(67,81,112)(68,180,129)(70,120,84)(71,143,124)(72,149,96)(75,187,222)(78,156,177)\\
(83,179,240)(85,181,100)(87,95,111)(88,214,159)(89,163,211)(93,183,133)(98,243,190)\\
(101,184,160)(102,220,130)(103,154,223)(104,202,150)(105,146,114)(106,178,194)\\
(113,155,132)(115,196,127)(118,171,123)(126,231,157)(135,226,200)(139,199,174)\\
(140,203,204)(148,207,192)(197,210,232)(215,227,241)(218,233,239)(244,473,409)\\
(245,339,388)(246,423,449)(247,430,471)(248,253,474)(249,254,470)(250,321,269)\\
(251,479,443)(252,311,419)(255,309,460)(256,279,371)(257,268,446)(258,438,264)\\
(259,389,363)(260,399,323)(261,319,336)(262,462,271)(263,295,267)(265,458,348)\\
(266,418,296)(270,393,333)(272,408,476)(273,310,324)(274,381,445)(275,297,379)\\
(276,404,289)(277,350,444)(278,451,298)(280,453,390)(281,434,299)(282,328,396)\\
(283,357,386)(284,288,358)(285,420,307)(286,417,454)(287,374,367)(290,485,313)\\
(291,373,326)(292,372,415)(293,304,340)(294,377,361)(300,480,335)(301,318,353)\\
(302,459,370)(303,375,481)(305,447,346)(306,400,337)(308,380,368)(312,468,484)\\
(314,402,369)(315,316,397)(317,410,391)(320,349,436)(322,352,366)(325,394,477)\\
(327,478,411)(329,398,401)(330,437,356)(331,359,413)(332,429,466)(334,482,387)\\
(338,464,362)(341,422,344)(342,457,354)(343,469,452)(345,463,385)(347,472,467)\\
(351,364,392)(355,455,450)(360,442,428)(365,395,456)(376,483,486)(378,412,431)\\
(382,448,383)(384,403,406)(405,425,426)(407,424,414)(416,433,427)(421,461,465)\\
(432,441,475)(435,440,439) }

\section*{Acknowledgements}

The authors acknowledge financial support from an NSERC Discovery Grant, a Memorial University of Newfoundland startup grant, and a grant from the Memorial University of Newfoundland Seed, Bridge and Multidisciplinary Fund (all held by the first author).  This work has been supported in part by the Croatian Science Foundation under the project 6732.  The first author thanks the Department of Mathematics at the University of Rijeka for their hospitality.


\begin{thebibliography}{99}

\bibitem{incidence} R.~F.~Bailey, On the metric dimension of incidence graphs, {\em Discrete Math.} {\bf 341} (2018), 1613--1619.

\bibitem{distanceregular.org} R.~F.~Bailey, F.~L.~Bosquet, C.~M.~Bowers, A.~D.~M.~Jackson and C.~H.~Weir, \url{https://www.distanceregular.org}, 2017--present.

\bibitem{BvLS73} E.~R.~Berlekamp, J.~H.~van Lint and J.~J.~Seidel, A strongly regular graph derived from the perfect ternary Golay code, in {\em A Survey of Combinatorial Theory (Proc. Internat. Sympos., Colorado State Univ., Fort Collins, Colo., 1971)}, pp. 25--30. North-Holland, Amsterdam, 1973. 

\bibitem{BJL} T.~Beth, D.~Jungnickel and H.~Lenz, {\em Design Theory} (second edition), Volume I, Cambridge University Press, Cambridge, 1999.

\bibitem{Magma} W.~Bosma, J.~J.~Cannon and C.~Playoust, The Magma algebra system. I. The user language, {\em J. Symbolic Comput.} {\bf 24} (1997), 235--265.

\bibitem{BCN} A.~E.~Brouwer, A.~M.~Cohen and A.~Neumaier, {\em Distance-Regular Graphs}, Springer-Verlag, Berlin, 1989.

\bibitem{BCN-corr-add} A.~E.~Brouwer, A.~M.~Cohen and A.~Neumaier, Corrections and additions to the book `Distance-Regular Graphs', available from \url{https://homepages.cwi.nl/~aeb/math/bcn/}.

\bibitem{BH92} A.~E.~Brouwer and W.~H.~Haemers, Structure and uniqueness of the $(81,20,1,6)$ strongly regular graph, {\em Discrete Math.} {\bf 106/107} (1992), 77--82.

\bibitem{BKR98} A.~E.~Brouwer, J.~H.~Koolen and R.~J.~Riebeek, A new distance-regular graph associated to the Mathieu group $M_{10}$, {\em J. Algebraic Combin.} {\bf 8} (1998), 153--156.

\bibitem{BvL84} A.~E.~Brouwer and J.~H.~van Lint, Strongly regular graphs and partial geometries, in {\em Enumeration and Design (Waterloo, Ont., 1982)}, eds. D.~M.~Jackson and S.~A.~Vanstone, pp. 85--122. Academic Press, Toronto, 1984. 

\bibitem{vDKT} E.~R.~van Dam, J.~H.~Koolen and H.~Tanaka, Distance-regular graphs, {\em Electronic J. Combin.} (2016), Dynamic Survey \#DS22.

\bibitem{GAP} The GAP~Group, \emph{GAP -- Groups, Algorithms, and Programming, Version 4.10.2} (2019); \url{https://www.gap-system.org}.

\bibitem{ILPP97} A.~A.~Ivanov, R.~A.~Liebler, T.~Penttila and C.~E.~Praeger, Antipodal distance-transitive covers of complete bipartite graphs, {\em European J. Combin.} {\bf 18} (1997), 11--33.

\bibitem{Jungnickel84} D.~Jungnickel, The number of designs with classical parameters grows exponentially, {\em Geom. Dedicata} {\bf 16} (1984), 167--178.

\bibitem{MacWS77} F.~J.~MacWilliams and N.~J.~A.~Sloane, {\em The Theory of Error-Correcting Codes}, North Holland, Amsterdam, 1977.

\bibitem{Mavron2000} V.~C.~Mavron and V.~D.~Tonchev, On symmetric nets and generalized Hadamard matrices from affine designs, {\em J. Geom.} {\bf 67} (2000), 180--187.

\bibitem{PS97} C.~E.~Praeger and L.~H.~Soicher, {\em Low rank representations and graphs for sporadic groups}, Australian Mathematical Society Lecture Series {\bf 8}, Cambridge University Press, Cambridge, 1997.

\bibitem{Soicher93} L.~H.~Soicher, Three new distance-regular graphs, {\em European J. Combin.} {\bf 14} (1993), 501--505.

\bibitem{DESIGN} L.~H.~Soicher, \emph{DESIGN: The Design Package for GAP, Version 1.7} (2019);\linebreak \url{https://gap-packages.github.io/design}.

\bibitem{GRAPE} L.~H.~Soicher, \emph{GRAPE: GRaph Algorithms using PErmutation groups, Version 4.8.2} (2019); \url{https://gap-packages.github.io/grape}.

\end{thebibliography}
\end{document}